\documentclass[12pt,twoside]{article}

\date{}

\begin{document}

\centerline {\Large{\bf Universality properties of double series  }}

\centerline{}

\centerline{\Large{\bf by generalized Walsh system }}

\centerline{}

\centerline{\bf {Sergo A. Episkoposian}}

\centerline{}

\centerline{Faculty of Applied Mathematics}

\centerline{State Engeniering University of Armenia}

\centerline{Yerevan, Teryan st.105,  375049, Armenia.}

\centerline{}

\newtheorem{Theorem}{\quad Theorem}[section]

\newtheorem{Definition}[Theorem]{\quad Definition}

\newtheorem{Corollary}[Theorem]{\quad Corollary}

\newtheorem{Lemma}[Theorem]{\quad Lemma}

\newtheorem{Property}[Theorem]{\quad Property}

\newtheorem{Example}[Theorem]{\quad Example}

\begin{abstract}
In this paper we consider a question on existence of double  
series by generalized Walsh system, which are  universal in weighted $L_\mu^1[0,1]^2$ spaces. In particular, we construct a weighted function $\mu(x,y)$ and a  double series by generalized Walsh system of the form
$$\sum_{n,k=1}^\infty c_{n,k}\psi_n(x)\psi_k(y)\ \  \mbox{with} \ \ \sum_{n,k=1}^\infty \left |
c_{n,k} \right|^q <\infty\ \mbox{for all}\ q>2,$$ which is universal
in $L_\mu^1[0,1]^2$ concerning subseries with respect to
convergence, in the sense of both spherical and rectangular partial
sums.

\end{abstract}



\section{Introduction}

Let $X$ be a Banach space.
\begin{Definition}
 A series
\begin{equation}
\sum_{k=1}^\infty f_k,\ \ f_k \in X
\end{equation}
is said to be universal in $X$ with respect to rearrangements, if
for any $f \in X$ the members of (1) can be rearranged so that
the obtained series $\displaystyle \sum_{k=1}^\infty
f_{\sigma(k)}$ converges to $f$ by norm of $X$.
\end{Definition}

\begin{Definition}
The series (1) is said to be universal (in $X$) concerning
subseries, if for any $f \in X$ it is possible to choose a
subseries $\displaystyle {\sum_{k=1}^\infty f_{n_k}}$ from
(1), which converges to the $f$ by norm of $X$.
\end{Definition}

 Note, that for one-dimensional case there are many papers
are devoted to the question on existence of various types of
universal series in the sense of convergence almost everywhere and
on a measure ( see [1] - [9], [11]).

The first usual universal in the sense of convergence almost
everywhere trigonometric series were constructed by D.E.Menshov
[1] and V.Ya.Kozlov [2]. The series of the form
\[
{1\over2}+\sum_{k=1}^\infty a_k \cos{kx}+b_k \sin{kx} 
\]
was constructed just by them such that for any measurable on
$[0,2\pi]$  function $f(x)$ there exists the growing sequence of
natural numbers $n_k$  such that the series  having the
sequence of partial sums with numbers $n_k$   converges to $f(x)$
almost everywhere on $[0,2\pi]$.

Note that in this result, when $f(x)\in{L^1_{[0,2\pi]}} $, it is
impossible to replace convergence almost everywhere by convergence
in the metric ${L^1_{[0,2\pi]}}$.

 This result was distributed by A.A.Talalian on arbitrary orthonormal
 complete systems . He also established (see [3]), that if
 $\{\phi_n(x)\}_{n=1}^\infty $  - the normalized basis of space ${L^p_{[0,1]}},p>1 $,
  then there exists a series of the form
\begin{equation} 
\sum_{k=1}^\infty{a_k\phi_k(x)},\ \ a_k \to 0.
\end{equation}
which has property: for any measurable function $f(x)$ the members
of series (2) can be rearranged so that the again received series
converge on a measure on [0,1] to $f(x)$.

In [4] O. P. Dzangadze these results are transferred to two-dimensional case.

W. Orlicz [5] observed the fact that there exist functional series
that are universal with respect to rearrangements in the sense of
a.e. convergence in the class of a.e. finite measurable functions.
It is also useful to note that even Rieman proved that every
convergent numerical series which is not absolutely convergent is
universal with respect to rearrangements in the class of all real
numbers.

Let $\mu(x)$, $0<\mu(x) \le1, x\in[0,1]$
 be a measurable on $[0,1]$ function and let $L_\mu^1[0,1]$ be a space of real measurable functions
 $f(x)$, $x\in [0,1]$ with
\[
  \int_0^1 |f(x)| \mu(x) dx<\infty.
  \]

In [6] - [9] it is proved the existence of universal one-dimensional series by trigonometric and classical Walsh system with respect to rearrangements and subseries. Some results for two-dimensional case for classical Walsh system was obtained in [11]. In this paper we consider this problems for double series by generalized Walsh system

\section{Preliminary Notes}

Let $a$ denote a fixed integer, $a\geq 2$ and put $\omega_a=e^{2
\pi i \over a}$.

Now we will give the definitions of generalized Rademacher and
Walsh systems (see [12] ).

\begin{Definition}
 The Rademacher  system of order $a$ is defined
by
$$\varphi_0(x)=\omega_a^k\ \ if \ \ x \in \left[ {k\over a}, {k+1\over a}\right),\ \ k=0,1,...,a-1,\ \ x\in[0,1)$$
and for $n\geq 0$
$$\varphi_n(x+1)=\varphi_n(x)=\varphi_0(a^nx).$$
\end{Definition}

\begin{Definition}
 The generalized Walsh system of order $a$ is
defined by
$$\psi_0(x)=1,$$
and if $n=\alpha_1a^{n_1}+...+\alpha_s a^{n_s}$ where
$n_1>...>n_s,$ then
$$\psi_n(x)=\varphi_{n_1}^{\alpha_1}(x)\cdot...\cdot \varphi_{n_s}^{\alpha_s}(x) .$$
\end{Definition}
 Let's denote the generalized Walsh system of order $a$ by $\Psi_{a}$, $a\geq 2$.

Note that $\Psi_2$ is the classical Walsh system.

The basic properties of the generalized Walsh system of order $a$
are obtained by   H.E.Chrestenson, R. Pely, J. Fine, W. Young, C.
Vatari, N. Vilenkin and others (see [12]- [17]).

First we present some properties of $\Psi_{a}$ system (see
Definition 2.1).

{\bf Property 1.} Each  $n$th Rademacher function  has period
$1\over a^n$ and
\begin{equation}
\varphi_n(x)=const \in \Omega_a=\{ 1, \omega_a,\omega_a^2,....,
\omega_a^{a-1} \},
\end{equation}
if  $ x\in \Delta_{n+1}^{(k)}= \left[ {k\over a^{n+1}}, {k+1\over
a^{n+1}}\right)$, $ k=0,...,a^{n+1}-1$, $ n=1,2,....$.

It is also easily verified, that
\begin{equation}
\left( \varphi_n(x) \right)^k =\left( \varphi_n(x) \right)^m, \ \
\forall n,k\in \mathcal{N},  \textit{where}\ \  m=k\ (\textit{mod}
\ a)
\end{equation}

{\bf Property 2.} It is clear, that for any integer $n$ the Walsh
function $\psi_n(x)$ consists of a finite product of Rademacher
functions and accepts values from $\Omega_a$.

{\bf Property 3.} The generalized Walsh system $\Psi_a$, $a\geq 2$
is a complete orthonormal system in $L^2[0,1)$ and basis in
$L^p[0,1]$, $p>1$ (see [5]).

The rectangular and spherical partial sums of the double series
\[
\sum_{k,\nu=1}^\infty c_{k,\nu} \psi_k(x)\psi_\nu(y)
\]
will be denoted by
\[
 S_{n,m}(x,y)= \sum _{k=1}^n\sum _{\nu=1}^m c_{k,\nu} \psi_k(x)\psi_\nu(y)
\]
and
\[
 S_R(x,y)= \sum _{\nu^2+k^2 \leq R^2} c_{k,\nu} \psi_k(x)\psi_\nu(y).
\]

If $g(x,y)$ is a continuous function on $T=[0,1]^2$, then we set
 \[
 ||g(x,y)||_C=\max_{(x,y)\in T} |g(x,y)|.
 \]
\section{Main Results}

Let's denote the generalized Walsh system of order $a$ by
$\Psi_{a}$, $a\geq 2$.
 These are the main results of the paper.

\begin{Theorem}
 There exists a double series of the form
\begin{equation}
\sum_{n,k=1}^\infty{c_{n,k}\psi_n(x)\psi_k(y) \ \  with \ \
\sum_{n,k=1}^\infty \left | {c_{n,k}} \right|^q <\infty} \ \ for \
all\ q>2
\end{equation}
with the following property: \noindent for any number
$\varepsilon>0$ a weighted function $\mu(x,y)$ satisfying
\begin{equation}
0<\mu(x,y) \le 1, \left | \{ (x,y)\in T: \mu(x,y)\not =1
\} \right | <\varepsilon
\end{equation}
  can be constructed so that the series (5) is universal
  in $L_\mu^1(T)$ concerning subseries
with respect to convergence in the sense of both spherical and
rectangular partial sums.
 \end{Theorem}

\begin{Theorem}  There exists a double series of the form
(5) with the following property: \noindent for any number
$\varepsilon>0$ a weighted function $\mu(x,y)$ with (6) can be
constructed, so that the series (5) is universal in $L_\mu^1(T)$
concerning rearrangements  with respect to convergence in the sense
of both spherical and rectangular partial sums.
\end{Theorem}

Repeating the reasoning of the proof of Lemma 2 in [10] we'll
receive the following lemma:

\begin{Lemma} For any given numbers $0<\varepsilon<1$, $N_0>2$ $(N_0\in
\mathcal{N})$ and a step function
$$f(x)= \sum_{s=1}^q \gamma_s \cdot \chi_{\Delta_s} (x),$$
where  $\Delta_s$ is an interval  of the form $\displaystyle
\Delta_m^{(i)}= \left[ {{i-1}\over {2^m}},{i\over {2^m}} \right]
$, $ 1\leq i \leq 2^m $, there exist a measurable set $E \subset
[0,1]$ and a polynomial $P(x)$ of the form
 $$P(x)= \sum_{k=N_0}^N c_k\psi_k(x) $$
which satisfy the conditions:
$$P(x)=f(x)\ \  on \ \ E, \leqno(1)$$
$$|E|> (1- \varepsilon), \leqno(2)$$
$$\sum_{k=N_0}^N |c_k|^{2+ \varepsilon}< \varepsilon, \leqno(3)$$
$$\max_{N_0 \leq m<N} \left[ \int_e \left | \sum_{k=N_0}^m c_k
\psi_k(x) \right | dx \right] <\varepsilon+\int_e
|f_(x)|dx,\leqno(4)$$ for every measurable subset $e$  of $E$.
\end{Lemma}

Then applying this Lemma we get next one:

\begin{Lemma}
 For any numbers $\gamma \not=0$, $0<\delta<1$,
$N>1$ and for any square $\Delta=\Delta_1 \times \Delta_2 \subset
T$ there exists a measurable set $E \subset T$ and a
polynomial $P(x,y)$ of the form
$$P(x,y)= \sum_{k,s=N}^M c_{k,s}\psi_k(x) \cdot \psi_s(y) ,$$
with the following properties:
$$|E|> 1- \delta, \leqno(1)$$
$$\sum_{k,s=N}^M |c_{k,s}|^{2+\delta}< \delta, \leqno(2)$$
$$P(x,y)=\gamma \cdot \chi_\Delta (x,y) \ \  for \ \ (x,y) \in E, \leqno(3)$$
$$\max_{N \leq \overline {n}, \overline{m} \leq M} \left[\int \int_e \left | \sum_{k,s=N}^{\overline {n}, \overline{m}}  c_{k,s}\psi_k(x) \cdot \psi_s(y) \right | dx dy \right]\leqno(4)$$
$$+\max_{\sqrt {2}N \leq R \leq \sqrt {2}M} \left[ \int\int_e \left | \sum_{2N^2 \leq k^2+s^2 \leq R^2}  c_{k,s}\psi_k(x) \cdot \psi_s(y) \right | dx dy \right] \leq 16 \cdot |\gamma|\cdot |\Delta|,$$
for every measurable subset $e$  of $E$.
\end{Lemma}

{\bf Proof .} We apply Lemma 3.3, setting
$$f(x)= \gamma \cdot \chi_{\Delta_1} (x),\ \ N_0=N,\ \ \varepsilon= {\delta \over 2}.$$
Then we can define a measurable set $E_1 \subset [0,1]$ and a
polynomial $P_1(x)$ of the form
 $$P_1(x)= \sum_{k=N}^{N_1} a_k\psi_k(x) $$
which satisfy the conditions:
$$P_1(x)=\gamma \cdot \chi_{\Delta_1} (x) \ \  for \ \ x \in E_1, \leqno(1^0)$$
$$|E_1|> 1- {\delta \over 2}, \leqno(2^0)$$
$$\sum_{k=N}^{N_1} |a_k|^{2+\delta}< \delta, \leqno(3^0)$$
$$\max_{N \leq \overline {n} \leq N_1} \left[ \int_{e_1} \left | \sum_{k=N}^{\overline {n}}  a_k\psi_k(x)  \right | dx \right] \leq 2 \cdot |\gamma| \cdot |\Delta_1|, \leqno(4^0)$$
for every measurable subset $e_1$  of $E_1$.

Set
\begin{equation}
M_0=2 \cdot \left( N_1^2+1 \right)
\end{equation}
and apply Lemma 3.3 again, setting
$$f(y)=  \chi_{\Delta_2} (y),\ \ N_0=M_0,\ \ \varepsilon= {\delta \over 2}.$$
Then we can define a measurable set $E_2 \subset [0,1]$ and a
polynomial $P_2(y)$ of the form
 $$P_2(y)= \sum_{s=M_0}^M b_s\psi_s(y), $$
which satisfy the conditions:
$$P_2(y)= \chi_{\Delta_2} (y) \ \  for \ \ y \in E_2, \leqno(1^{00})$$
$$|E_2|> 1- {\delta \over 2}, \leqno(2^{00})$$
$$\sum_{s=M_0}^M |b_s|^{2+\delta}< \delta, \leqno(3^{00})$$
$$\max_{M_0 \leq \overline {m} \leq M} \left[ \int_{e_2} \left | \sum_{s=M_0}^{\overline {m}}  b_s\psi_s(y)  \right | dy \right] \leq 2 \cdot  |\Delta_2|, \leqno(4^{00})$$
for every measurable subset $e_2$  of $E_2$.

Set
\begin{equation}
E=E_1 \times E_2,
\end{equation}
\begin{equation}
P(x,y)= P_1(x) \cdot P_2(x)=\sum_{k,s=N}^M c_{k,s}\psi_k(x) \cdot
\psi_s(y) ,
\end{equation}
where
\begin{equation}
c_{k,s}=a_k \cdot b_s,\ \ if \ \ N \leq k \leq N_1, \ \ M_0 \leq s
\leq M
\end{equation}
and
\[
  c_{k,s}=0,\ \ for \ \ other \ \ k\ \ and\ \  s.
\]
By $(1^0) - (3^0),\ \ (1^{00}) - (3^{00})$ and (3.2) - (3.4) we
obtain
$$|E|> 1- \delta, $$
$$\sum_{k,s=N}^M |c_{k,s}|^{2+\delta}= \sum_{k=N}^{N_1}
|a_k|^{2+\delta} \cdot \sum_{s=M_0}^M |b_s|^{2+\delta}< \delta, $$
$$P(x,y)=\gamma \cdot \chi_\Delta (x,y) \ \  for \ \ (x,y) \in E.$$
Thus, the statements 1) - 3) of Lemma 3.4 are satisfied. Now we will
check the fulfillment of statement 4).

Let $N^2+M_0^2<R^2<N_1^2+M^2$, then for some $m_0>M_0$ we have
$m_0<R^2<m_0+1$ and from (3.1) it follows, that
$R^2-N_1^2>(m_0-1)^2$.

Consequently taking relations $(4^0), (4^{00})$ and (3.2) - (3.4)
for any measurable set $e \subset E$ $\left( e=e_1  \times e_2,\ \
e_1 \subset E_1,\ \ e_2 \subset E_2 \right)$ we obtain
$$ \int\int_e \left | \sum_{N^2+M^2 \leq
k^2+s^2 \leq R^2}  c_{k,s}\psi_k(x) \cdot \psi_s(y) \right | dx dy $$
$$\leq \int\int_e \left |
\sum_{k=N}^{N_1}\sum_{s=M_0}^{m_0-1} c_{k,s}\psi_k(x) \cdot \psi_s(y)
\right | dx dy $$
$$+\max_{N<n \leq
N_1}\left[ \int\int_e \left | \sum_{k=N}^n  c_{k,m_0}\psi_k(x) \cdot
\psi_{m_0}(y) \right | dx dy\right] $$
$$\leq \left[ \int_{e_1}
\left | \sum_{k=N}^{N_1} a_k\psi_k(x) \right | dx \right] \cdot \left[
\int_{e_2} \left | \sum_{s=M_0}^{m_0-1}  b_s\psi_s(y) \right | dy
\right]$$
$$+|b_{m_0}| \cdot  \left[ \int_{e_2} \left |
\psi_{m_0}(y) \right | dy \right] \cdot \max_{N<n \leq  N_1}\left[
\int_{e_1} \left | \sum_{k=N}^n a_k\psi_k(x) \right | dx \right] $$
$$\leq 12 \cdot |\gamma|\cdot |\Delta|.$$

Similarly, for $N \leq \overline {n} \leq N_1, \ \ M_0 \leq
\overline {m} \leq M$, we get
$$\int\int_e \left |
\sum_{k,s=N}^{\overline {n}, \overline{m}}  c_{k,s}\psi_k(x) \cdot
\psi_s(y) \right | dx dy  \leq 4\cdot |\gamma|\cdot |\Delta|.$$

{\bf Lemma 3.4 is proved.}
\par\par\bigskip

\begin{Lemma}
 For any numbers $\varepsilon>0$, $N>1$ and a step
function
$$  {f(x,y)= \sum_{\nu=1}^{\nu_0} \gamma_\nu \cdot
\chi_{\Delta_\nu} (x,y)},$$ there exists a measurable set $E \subset
T$ and a polynomial $P(x,y)$ of the form
\[
P(x,y)= \sum_{k,s=N}^M c_{k,s}\psi_k(x) \cdot \psi_s(y) ,
\]
which satisfy the following conditions:
$$ P(x,y)=f(x,y) \ \  for \ \ (x,y) \in E,\leqno(1^0)$$
$$ |E|> 1- \varepsilon, \leqno(2^0)$$
$$\sum_{k,s=N}^M |c_{k,s}|^{2+\varepsilon}< \varepsilon,\leqno(3^0)$$
$$ \max_{N \leq \overline {n}, \overline{m} < M} \left[ \int\int_e
\left | \sum_{k,s=N}^{\overline {n}, \overline{m}}  c_{k,s}\psi_k(x)
\cdot \psi_s(y) \right | dx dy \right]\leqno(4^0)$$
\[
+\max_{\sqrt {2}N \leq R \leq \sqrt {2}M} \left[ \int\int_e \left |
\sum_{2N^2 \leq k^2+s^2 \leq R^2}  c_{k,s}\psi_k(x) \cdot \psi_s(y) \right
| dx dy \right]
\]
\[
\leq 2 \cdot \int\int_e |f(x,y)| dxdy+ \varepsilon,
 \]
 for every measurable subset $e$  of $E$.
\end{Lemma}

{\bf Proof .} Without any loss of generality, we assume
that
\begin{equation}
\max_{1 \leq \nu \leq \nu_0} \left( |\gamma_\nu| \cdot |\Delta_\nu|
\right)< {\varepsilon \over 32},
\end{equation}
($\Delta_ \nu,\ \ 1 \leq \nu \leq \nu_0$ are the constancy
rectangular domian of $f(x,y)$, i.e. where the function $f(x,y)$ is
constant).

Given an integer $1 \leq \nu \leq \nu_0$, by applying Lemma 3.4 with
$\delta={\varepsilon \over 16 \nu_0}$, we find that there exists a
measurable set $E_\nu \subset T$ and a polynomial $P_\nu(x,y)$ of
the form
\begin{equation}
P_\nu(x,y)= \sum_{k,s=N_\nu}^{M_\nu} c_{k,s}^{(\nu)}\psi_k(x) \cdot
\psi_s(y)
\end{equation}
with the following properties:
\begin{equation}
|E_\nu |> 1- {\varepsilon \over 2^\nu},
\end{equation}
\begin{equation}
\sum_{k,s=N_\nu}^{M_\nu} |c_{k,s}^{(\nu)}|^{2+\varepsilon}<
{\varepsilon \over \nu_0},
\end{equation}
\begin{equation}
P_\nu (x,y)=\gamma_\nu \cdot \chi_{\Delta_\nu} (x,y) \ \  for \ \
(x,y) \in E_\nu,
\end{equation}
$$\max_{N_\nu \leq \overline {n}, \overline{m} \leq M_\nu} \left[
\int\int_e \left | \sum_{k,s=N_\nu}^{\overline {n}, \overline{m}}
c_{k,s}^{(\nu)}\psi_k(x) \cdot \psi_s(y) \right | dx dy \right]$$
$$+\max_{\sqrt {2}N_\nu \leq R \leq \sqrt {2}M_\nu} \left[
\int\int_e \left | \sum_{2N_\nu^2 \leq k^2+s^2 \leq R^2}
c_{k,s}^{(\nu)}\psi_k(x) \cdot \psi_s(y) \right | dx dy \right] $$
\begin{equation}
\leq 16 \cdot |\gamma_\nu|\cdot |\Delta_\nu|<{\varepsilon \over 2},
\end{equation}
for every measurable subset $e$ of $E_\nu$ (see (11)).

Then we can take
$$N_1=N,\ \ n_\nu=M_{\nu-1}+1,\ \ 1 \leq \nu \leq \nu_0.$$
Set
\begin{equation}
E=\bigcap_{\nu=1}^{\nu_0} E_\nu,
\end{equation}
\begin{equation}
P(x,y)= \sum_{\nu=1}^{\nu_0}P_\nu(x,y)=\sum_{k,s=N}^M c_{k,s}\psi_k(x)
\cdot \psi_s(y),\ \ M=M_{\nu_0},
\end{equation}
where
\begin{equation}
c_{k,s}=c_{k,s}^{(\nu)},\ \ for \ \ N_\nu \leq k,s \leq M_\nu, \ \ 1
\leq \nu \leq \nu_0
\end{equation}
and
\[
 c_{k,s}=0,\ \ for \ \ other \ \ k\ \ and\ \ s.
\]
From (13) - (15), (17) - (19) we obtain:
$$P(x,y)=f(x,y) \ \  for \ \ (x,y) \in E, $$
$$|E|> 1- \varepsilon, $$
$$\sum_{k,s=N}^M |c_{k,s}|^{2+\varepsilon}< \sum_{\nu=1}^{\nu_0}
\left[ \sum_{k,s=N_\nu}^{M_\nu} |c_{k,s}^{(\nu)}|^{2+\varepsilon}
\right]< \varepsilon $$.

Then, let $R \in [\sqrt {2}N,\sqrt {2}M ]$, then for some $\nu', 1 \leq
\nu' \leq \nu_0$ we have $\sqrt {2}N_{\nu'} \leq R \leq \sqrt
{2}N_{\nu'+1}$, consequently from (18) and (19) we have
\[
\sum_{2N^2 \leq k^2+s^2 \leq R^2}  c_{k,s}\psi_k(x) \cdot
\psi_s(y)=\sum_{\nu=1}^{\nu'-1}P_\nu(x,y)
 \]
 \[
 +\sum_{2N_{\nu'}^2 \leq
k^2+s^2 \leq R^2}  c_{k,s}^{(\nu')}\psi_k(x) \cdot \psi_s(y).
\]

In view of the conditions (13) - (16) and the equality
$P(x,y)=f(x,y)$ on $E$, for any measurable set $e \subset E$ we
obtain
\[
 \int\int_e \left | \sum_{2N^2 \leq k^2+s^2 \leq R^2}
c_{k,s}\psi_k(x) \cdot \psi_s(y) \right | dxdy
\]
\[
\leq \int\int_e \left | \sum_{\nu=1}^{\nu'-1}P_\nu(x,y) \right |
dxdy
\]
\[
+ \int\int_e \left | \sum_{2N_{\nu'}^2 \leq k^2+s^2 \leq R^2}
c_{k,s}^{(\nu')}\psi_k(x) \cdot \psi_s(y) \right | dxdy
\]
\[
\leq \int\int_e |f(x,y)|dxdy+{\varepsilon \over 2}.
\]
Similarly, for any $e \subset E$ we have
\[
\max_{N \leq \overline {n}, \overline{m} \leq M} \left[ \int\int_e
\left | \sum_{k,s=N}^{\overline {n}, \overline{m}}  c_{k,s}\psi_k(x)
\cdot \psi_s(y) \right | dx dy \right]
\]
\[
 \leq \int\int_e
|f(x,y)|dxdy+{\varepsilon \over 2}.
\]

{\bf Lemma 3.5 is proved.}
\par\par\bigskip

\section{Proofs of the theorems}

\vskip 4mm
The Theorem 3.1 is proved similarly Theorem 3 in [11], but for maintenance of integrity of this paper, here we will give the proof :
\vskip 4mm
{\bf Proof of Theorem 3.1.}
 \vskip 4mm Let
\begin{equation}
\{ f_s(x,y)\}_{s=1}^\infty,\ \  (x,y) \in T
\end{equation}
 be a sequence of all step functions, values and constancy interval endpoints of which
are rational numbers. Applying Lemma 3.5 consecutively, we can find a
sequence
  $\{ E_s\}_{s=1}^\infty $ of sets and a sequence of polynomials
\begin{equation}
P_s(x,y)=\sum_{k,\nu=N_{s-1}}^{ N_s-1}c_{k,\nu}^{(s)}\psi_k(x)\psi_\nu(y)
,
\end{equation}
\[
 1=N_0<N_1<...<N_s<....,\ \ s=1,2,....,
\]
which satisfy the conditions:
\begin{equation}
P_s(x,y)=f_s(x,y),\ \ (x,y)\in E_s,
\end{equation}
\begin{equation}
\left| E_s\right| >1-2^{-2(s+1)} ,\ \  E_s\subset T,
\end{equation}
\begin{equation}
\sum_{k,\nu=N_{s-1}}^{ N_s-1}\left
|c_{k,\nu}^{(s)}\right|^{2+2^{-2s}}< 2^{-2s},
\end{equation}
\[
\max_{N_{s-1} \leq \overline {n}, \overline{m} < N_{s}} \left[
\int\int_e \left | \sum_{k,\nu=N_{s-1}}^{\overline {n},
\overline{m}} c_{k,\nu}^{(s)}\psi_k(x) \cdot \psi_\nu(y) \right | dx dy
\right]
\]
\[+\max_{\sqrt {2}N_{s-1} \leq R \leq \sqrt {2}N_s} \left[ \int\int_e \left
| \sum_{2N_{s-1}^2 \leq k^2+\nu^2 \leq R^2}  c_{k,\nu}^{(s)}\psi_k(x)
\cdot \psi_\nu(y) \right | dx dy \right]
\]
\begin{equation}
\leq 2 \cdot \int\int_e |f_s(x,y)| dxdy+ 2^{-2(s+1)},
\end{equation}
for every measurable subset $e$  of $E_s$.

Denote
\begin{equation}
\sum_{k,\nu=1}^\infty c_{k,\nu} \psi_k(x)\psi_\nu(y)=\sum_{s=1}^\infty
\left[\sum_{k,\nu=N_{s-1}}^{ N_s-1}c_{k,\nu}^{(s)}\psi_k(x)\psi_\nu(y)
\right],
\end{equation}
where
\[
c_{k,\nu}=c_{k,\nu}^{(s)}, \ \  for\ \ N_{s-1}\leq k,\nu <N_s,\ \
s=1,2,... .
\]

For an arbitrary number $\varepsilon >0$  we set
\begin{equation}
\cases { \Omega_n =\displaystyle{ \bigcap_{s=n}^\infty E_s,}\ \
n=1,2,.... ; \cr
 E=\Omega_{n_0} =\displaystyle{ \bigcap_{s=n_0}^\infty E_s,}\ \   n_0=[\log_{1/2} \varepsilon]+1;
\cr
 B= \displaystyle{\bigcup _{n=n_0} ^\infty} \Omega_n =\Omega_{n_0}
\bigcup \left( \displaystyle{\bigcup _{n=n_0+1}^ \infty }\Omega_n
\setminus \Omega_{n-1} \right). }
 \end{equation}
It is obvious ( see (23), (27) ) that $\left| B \right|=1$  and
$\left| E \right| >1- \varepsilon .$

We define a function $\mu(x,y)$ in the following way:
\begin{equation}
\mu(x,y)= \cases { 1, \ \ for \ \ (x,y) \in E \cup (T \setminus
B);\cr  \mu_n, \ \ for  \ \ (x,y) \in \Omega_n \setminus
\Omega_{n-1},\ \ n\geq n_0+1,}
\end{equation}
where
\begin{equation}
 \displaystyle{ \mu_n=\left[
2^{2n}\cdot \prod_{s=1}^n h_s \right]^{-1}};\hfil \quad
\end{equation}
\[
 h_s=|| f_s||_C+\displaystyle{ \max_{N_{s-1} \leq \overline {n}, \overline{m} < N_{s}}
 \left|\left|  \sum_{k,\nu=N_{s-1}}^{\overline {n},
\overline{m}} c_{k,\nu}^{(s)}\psi_k(x) \cdot \psi_\nu(y)\right|
\right|_C} \hfil
\]
\[
+ \displaystyle{ \max_{\sqrt {2}N_{s-1} \leq R \leq \sqrt {2}N_s}
\left|\left| \sum_{2N_{s-1}^2 \leq k^2+\nu^2 \leq R^2}
c_{k,\nu}^{(s)}\psi_k(x) \cdot \psi_\nu(y)\right| \right|_C +1.}
\]
From (24), (26) - (29) we obtain\\
$(A) - 0<\mu(x,y) \le1, \ \ \mu(x,y)$ is a measurable function and
\[
\left | \{(x,y)\in T:\mu(x,y)\not =1\} \right|<\varepsilon.
\]
(B) -- $\displaystyle {\sum_{k,\nu=1}^\infty
\left|c_{k,\nu}\right|^q<\infty} \ \mbox{for all}\ q>2.$

Hence, obviously we have (see (24) and (26))
$$ \lim_{\min\{k,\nu\}\to \infty}c_{k,\nu}=0. $$
It follows  from (27) - (29) that for all $s \geq n_0$  and
$N_{s-1} \leq \overline {n}, \overline{m} < N_{s}$
\[
\int\int_{T \setminus \Omega_s} \left|
\sum_{k,\nu=N_{s-1}}^{\overline {n}, \overline{m}}
c_{k,\nu}^{(s)}\psi_k(x) \cdot \psi_\nu(y)\right| \mu(x,y) dxdy
\]
\[
=\sum_{n=s+1}^ \infty \left[\int\int_{\Omega_n \setminus
\Omega_{n-1}} \left| \sum_{k,\nu=N_{s-1}}^{\overline {n},
\overline{m}} c_{k,\nu}^{(s)}\psi_k(x) \cdot \psi_\nu(y)\right| \mu_n dxdy
\right]
\]
\begin{equation}
\leq \sum_{n=s+1}^ \infty2^{-2n} \left[\int\int_T \left|
\sum_{k,\nu=N_{s-1}}^{\overline {n}, \overline{m}}
c_{k,\nu}^{(s)}\psi_k(x) \cdot \psi_\nu(y)\right| h_s^{-1} dxdy \right]<{1
\over 3}2^{-2s}.
\end{equation}
Analogously for all $s \geq n_0$  and $\sqrt {2}N_{s-1} \leq R \leq
\sqrt {2}N_s$ we have
\begin{equation}
\int\int_{T \setminus \Omega_s} \left|\sum_{2N_{s-1}^2 \leq
k^2+\nu^2 \leq R^2} c_{k,\nu}^{(s)}\psi_k(x) \cdot \psi_\nu(y)\right|
\mu(x,y) dxdy<{1 \over 3}2^{-2s}.
\end{equation}
By (21), (27) - (29) for all  $s \geq n_0$ we have
\[
\int\int_T \left| P_s(x,y)-f_s(x,y) \right|\mu(x,y)dxdy
\]
\[
=\int\int_{\Omega_s} \left| P_s(x,y)-f_s(x,y) \right|\mu(x,y)dxdy
\]
\[+\int\int_{T \setminus \Omega_{s}} \left| P_s(x,y)-f_s(x,y)
\right|\mu(x,y)dxdy
\]
\[
=\sum_{n=s+1}^\infty \left[\int\int_{\Omega_n \setminus
\Omega_{n-1}} \left| P_s(x,y)-f_s(x,y)  \right| \mu_n dxdy\right]
\]
\[
\leq \sum_{n=s+1}^ \infty 2^{-2n}\left[ \int\int_T \left(\left|
f_s(x,y) \right| + \sum_{k,\nu=N_{s-1}}^ {N_s-1}
c_{k,\nu}^{(s)}\psi_k(x) \cdot \psi_\nu(y) \right) h_s^{-1}dxdy \right ]
\]
\begin{equation}
<{1 \over 3}2^{-2s}<2^{-2s}.
\end{equation}
By (25) and (27) - (30) for all
 $N_{s-1}\leq \overline {n}, \overline{m} < N_{s}$ and $s
\geq n_0+1$ we obtain
\[
\int\int_T  \left| \sum_{k,\nu=N_{s-1}}^{\overline {n},
\overline{m}} c_{k,\nu}^{(s)}\psi_k(x) \cdot \psi_\nu(y)\right| \mu(x,y)
dxdy
\]
\[
\int\int_{\Omega_s}  \left| \sum_{k,\nu=N_{s-1}}^{\overline {n},
\overline{m}} c_{k,\nu}^{(s)}\psi_k(x) \cdot \psi_\nu(y)\right| \mu(x,y)
dxdy
\]
\[
\int\int_{T \setminus \Omega_s}  \left|
\sum_{k,\nu=N_{s-1}}^{\overline {n}, \overline{m}}
c_{k,\nu}^{(s)}\psi_k(x) \cdot \psi_\nu(y)\right| \mu(x,y) dxdy
\]
\[
< \sum_{n=n_0+1}^ s \left[\int\int_{\Omega_n \setminus \Omega_{n-1}}
\left| \sum_{k,\nu=N_{s-1}}^{\overline {n}, \overline{m}}
c_{k,\nu}^{(s)}\psi_k(x) \cdot \psi_\nu(y)\right|\cdot \mu_ndxdy
\right]+{1 \over 3}2^{-2s}
\]
\[
<\sum_{n=n_0+1}^ s \left( 2^{-2(s+1)}+2\cdot\int\int_{\Omega_n
\setminus \Omega_{n-1}} |f_s(x,y)|dxdy \right)\cdot \mu_n +{1 \over
3}2^{-2s}
 \]

\[
=2^{-2(s+1)} \cdot \sum_{n=n_0+1}^ s \mu_n+ +\int \int_{\Omega_s}
|f_s(x,y)|\mu(x,y)dxdy +{1 \over 3}2^{-2s}
\]
\begin{equation}
<2\cdot\int \int_T |f_s(x,y)|\mu(x,y)dxdy +2^{-2s}.
\end{equation}
Analogously for all $s \geq n_0$  and $\sqrt {2}N_{s-1} \leq R \leq
\sqrt {2}N_s$ we have (see (31))
\[
\int\int_T \left|\sum_{2N_{s-1}^2 \leq k^2+\nu^2 \leq R^2}
c_{k,\nu}^{(s)}\psi_k(x) \cdot \psi_\nu(y)\right| \mu(x,y) dxdy
\]
\begin{equation}
<2\cdot\int \int_T |f_s(x,y)|\mu(x,y)dxdy +2^{-2s}.
\end{equation}
Now we'll show that the series (26) is universal in $L_\mu^1(T)$
concerning subseries with respect to convergence by both spherical
and rectangular partial sums.

Let $ f(x,y) \in L_{\mu}^1 (T)$ , i. e.
\[
\int\int_T |f(x,y)|\mu(x,y) dxdy<\infty.
\]

It is easy to see that we can choose a function $f_{n_1}(x,y)$ from
the sequence (20) such that
\begin{equation}
\int\int_T \left| f(x,y)- f_{n_1}(x,y) \right|\mu(x,y)dxdy<2^{-2},\
\ n_1
> n_0+1.
\end{equation}
 Hence, we have
\begin{equation}
\int\int_T \left| f_{n_1}(x,y) \right|\mu(x,y)dxdy<2^{-2}+\int\int_T
|f(x,y)|\mu(x,y)dxdy.
 \end{equation}
 From (33) and (35) we get
\[
\int\int_T \left| f(x,y)- P_{n_1}(x,y) \right|\mu(x,y)dxdy
\]
\[
\leq \int\int_T \left| f(x,y)- f_{n_1}(x,y) \right|\mu(x,y)dxdy
\]
\[
+\int\int_T \left| f_{n_1}(x,y)- P_{n_1}(x,y) \right|\mu(x,y)dxdy<
2\cdot 2^{-2}.
 \]
Assume that numbers $n_1<n_2<...<n_{q-1}$ are chosen in such a way
that the following condition is satisfied:
\begin{equation}
\int\int_T \left| f(x,y)- \sum_{s=1}^j P_{n_s}(x,y)
\right|\mu(x,y)dxdy<2\cdot 2^{-2j}, \ \ 1\leq j \leq q-1 .
\end{equation}
Now we choose a function $f_{n_q}(x,y)$ from the sequence (20) such
that
\[
\int\int_T \left| \left( f(x,y)- \sum_{s=1}^{q-1} P_{n_s}(x,y)
\right)-f_{n_q}(x,y)\right| \mu(x,y)dxdy
\]
\begin{equation}
 <2\cdot 2^{-2q},\ \ n_q>n_{q-1}.
\end{equation}
This with (37) imply
\begin{equation}
\int\int_T \left| f_{n_q}(x,y)\right| \mu(x,y)dxdy<
 2^{-2q}+2\cdot 2^{-2(q-1)}=9\cdot 2^{-2q}.
\end{equation}
Hence and from (21), (32) - (34) we obtain
\begin{equation}
\int\int_T \left| f_{n_q}(x,y)- P_{n_q}(x,y) \right|\mu(x,y)dxdy<
 2^{-2n_q},
\end{equation}
where
\[
P_{n_q}(x,y)=\sum_{k,\nu=N_{n_q-1}}^{
N_{n_q}-1}c_{k,\nu}^{(n_q)}\psi_k(x)\psi_\nu(y) ,
\]
\[
\max_{N_{n_q-1} \leq \overline {n}, \overline{m}<N_{n_q}} \left[
\int\int_T \left| \sum_{k,\nu=N_{n_q-1}}^{\overline {n},
\overline{m}} c_{k,\nu}^{(n_q)}\psi_k(x) \cdot \psi_\nu(y)\right| \mu(x,y)
dxdy \right ]
\]
\begin{equation}
<19 \cdot 2^{-2q}.
\end{equation}
Analogously we have
\[
\max_{\sqrt {2}N_{n_q-1} \leq R \leq \sqrt {2}N_{n_q}} \left[
\int\int_T \left|\sum_{2N_{n_q-1}^2 \leq k^2+\nu^2 \leq R^2}
c_{k,\nu}^{(n_q)}\psi_k(x) \cdot \psi_\nu(y)\right| \mu(x,y) dxdy\right]
\]
\[
<19 \cdot 2^{-2q}.
\]
 In quality subseries of Theorem we shall take
\begin{equation}
\sum_{q=1}^\infty P_{n_q}(x,y)=\sum_{q=1}^\infty
\left[\sum_{k,\nu=N_{n_q-1}}^{
N_{n_q}-1}c_{k,\nu}^{(n_q)}\psi_k(x)\psi_\nu(y) \right].
\end{equation}
From (38) and (39) we have
\[
\int\int_T \left| f(x,y)-\sum_{s=1}^q P_{n_s}(x,y)\right| \mu(x,y)
dxdy
\]
\[
\leq \int\int_T \left|\left( f(x,y)- \sum_{s=1}^{q-1} P_{n_s}(x,y)
\right)-f_{n_q}(x,y)\right| \mu(x,y)dxdy
\]
\begin{equation}
+\int\int_T \left| f_{n_q}(x,y)- P_{n_q}(x,y) \right|\mu(x,y)dxdy<
 2\cdot 2^{-2q}.
\end{equation}

Let $ \overline {n}$ and $\overline{m}$ be arbitrary natural
numbers. Then for some natural number $q$ we have
\[
 N_{n_q-1}\leq \min \{\overline {n}, \overline{m}\}<N_{n_q}.
\]
Taking into account (40) and (42) for rectangular partial sums
$S_{\overline {n},\overline {m}}(x,y)$ of (41) we get
\[
\int\int_T \left| S_{\overline {n},\overline {m}}(x,y)-f(x,y)\right|
\mu(x,y) dxdy
\]
\[\leq \int\int_T \left|
f(x,y)-\sum_{s=1}^q P_{n_s}(x,y)\right| \mu(x,y) dxdy
\]
\[
+\max_{N_{n_q-1} \leq \overline {n}, \overline{m}<N_{n_q}} \left[
\int\int_T \left| \sum_{k,\nu=N_{n_q-1}}^{\overline {n},
\overline{m}} c_{k,\nu}^{(n_q)}\psi_k(x) \cdot \psi_\nu(y)\right| \mu(x,y)
dxdy \right ]
\]
\begin{equation}
 <21 \cdot 2^{-2q}.
\end{equation}
Analogously for $\sqrt {2}N_{n_q-1} \leq R \leq \sqrt {2}N_{n_q}$ we
have
\begin{equation}
\int\int_T \left|S_R(x,y)-f(x,y)\right| \mu(x,y) dxdy<21 \cdot
2^{-2q},
\end{equation}
where $S_R(x,y)$ the spherical partial sums of (41).

From (44) and (45) we conclude that the series (26) is
universal in $L_\mu^1(T)$ concerning subseries with respect to
convergence by both spherical and rectangular partial sums (see
Definition 1.2).

\vskip 4mm {\bf Theorem 3.1 is proved.}

\vskip 4mm

 {\bf Remark.} We can show Theorem 3.2 by the method in the proof of Theorem 3.1.


\begin{thebibliography}{99}

\bibitem{1}{D.E.Menshov, On the partial summs of trigonometric series,
 \em Studia Math., }{\bf 20} (1947), 197-238 [in russian].

\bibitem{2}{V.Ya.Kozlov , On the complete systems of orthogonal functions, \em Mat. Sbornik, }{\bf 26} (1950), 351-364 [in russian].

\bibitem{3}{A.A.Talalian, On the universal series with respect to rearrangements,
 \em Izv. Akad. Nauk SSSR, }{\bf 24} (1960), 567-604 [in russian].

\bibitem{4}{O.P.Dzangadze, On the universal double series,
 \em Bull. Georgian Acad. Sci. }{\bf 34} (1964), 225 - 228 [in russian].

\bibitem{5}{W.Orlicz , \"{U}ber die unabhangig von der Anordnung fast
     \"{u}berall kniwergenten Reihen, \em Bull. Acad. Polonaise Sci., }{\bf 81} (1927), 117-125.

\bibitem{6}{M.G.Grigorian , On the representation of functions by orthogonal series in weighted $L^p$ spaces, \em Studia Math., }{\bf 134} (1999), 211 - 237.

 \bibitem{7}{S.A. Episkoposian , On the series by Walsh system universal in weighted $L_\mu^1[0,1]$ spaces ,  \em Izv. Nats. Akad. Nauk Armenii, Math., English trans. in: J. Contemp. Math. Anal., }{\bf 34} (1999), 25 - 40.

\bibitem{8}{M.G.Grigorian, S.A. Episkoposian, Representation of functions
in weighted spaces $L_\mu^1[0,1]$ by  trigonometric and Walsh
series, \em Anal. Math., }{\bf 27} (2001), 267 - 277.

\bibitem{9}{S.A. Episkoposian , On the existence of universal series by trigonometric system,  \em J. Func. Anal., }{\bf 230} (2006), 169 - 189.

\bibitem{10}{S.A. Episkoposian, $L^1$- convergence of greedy algorithm by generalized Walsh system,
 \em Banach  Journal of  Mathematical Analysis, }{\bf 6} (2012), 161-174.

\bibitem{11}{S.A. Episkoposian, Existence of double Walsh series universal in
weighted spaces,  \em International Journal of Modern Mathematics,
}{\bf 2} (2007), 231 - 247.

\bibitem{12}{H.E. Chrestenson, A class of generalized Walsh functions,  \em  Pacific J. Math., }{\bf 45} (1955), 17 - 31.

\bibitem{13}{R. E. A. C. Paley, A remarkable set of orthogonal functions,  \em  London Math. Soc., }{\bf 34} (1932), 241 - 279.

\bibitem{14}{J. Fine, The generalized Walsh functions,  \em  Trans. Amer.
Math. Soc., }{\bf 69} (1950), 66 - 77.

\bibitem{15}{ W. Young, Mean convergence of generalized Walsh - Fourier series,  \em  Trans. Amer. Math. Soc., }{\bf 218} (1976), 311 - 320.

\bibitem{16}{C. Watari, On generalized Walsh-Fourier series,  \em  Tohoku Math. J., }{\bf 10} (1958), 211 - 241.

\bibitem{17}{ N. Vilenkini, On a class of complete orthonormal systems,  \em  AMS Transl., }{\bf 28} (1963), 1 - 35.


\end{thebibliography}
\end{document}